\documentclass[12pt]{amsart}
\usepackage{amsfonts, amsmath, amstext, amssymb}

\usepackage{latexsym}

\usepackage[english]{babel}
\newtheorem{thm}{Theorem}[section]

\newtheorem{cor}[thm]{Corollary}

\newtheorem{example}[thm]{Example}

\newtheorem{definition}[thm]{Definition}

\newcommand{\Q}{\mathbb{Q}}
\newcommand{\R}{\mathbb{R}}

\newcommand{\ml}{Martin-L\"{o}f }

\newcommand{\reals}{2^\omega}

\newcommand{\strings}{2^{<\omega}}
\newcommand{\nat}{\in\omega}

\newcommand{\dom}{\mbox{dom}}

\newcommand{\restr}{\!\!\restriction\!\!}
\newcommand{\st}{\;|\;}
\newcommand{\conv}{\!\!\downarrow}
\newcommand{\dvrg}{\!\!\uparrow}

\newcommand{\MLR}{\mathsf{MLR}}

\newcommand{\degd}{\mathbf{d}}

\newcommand{\ce}{c.e.\ }

\newcommand{\comment}[1]{}

\def\s01{\ensuremath{\Sigma^0_1}}
\def\d02{\ensuremath{\Delta^0_2}}
\def\phi{\varphi}

\begin{document}

\title{Key developments in algorithmic randomness}

\author[Franklin]{Johanna N.Y.\ Franklin}

\author[Porter]{Christopher P.\ Porter}

\date{June 10,2019}

\maketitle
\tableofcontents

\section{Introduction}\label{sec:intro}

The goal of this introductory survey is to present the major developments of algorithmic randomness with an eye toward its historical development. While two highly comprehensive books \cite{dhbook,niesbook} and one thorough survey article \cite{dhnt} have been written on the subject, our goal is to provide an introduction to algorithmic randomness that will be both useful for newcomers who want to develop a sense of the field quickly and interesting for researchers already in the field who would like to see these results presented in chronological order.

We begin in this section with a brief introduction to computability theory as well as the underlying mathematical concepts that we will later draw upon. Once these basic ideas have been presented, 
we will selectively survey four broad periods in which the primary developments in algorithmic randomness occurred:  (1) the mid-1960s to mid-1970s, in which the main definitions of algorithmic randomness were laid out and the basic properties of random sequences were established, (2) the 1980s through the 1990s, which featured intermittent and important work from a handful of researchers, (3) the 2000s, during which there was an explosion of results as the discipline matured into a full-fledged subbranch of computability theory, and (4) the early 2010s, which we briefly discuss as a lead-in to the remaining surveys in this volume, which cover in detail many of the exciting developments in this later period.

We do not intend this to be a full reconstruction of the history of algorithmic randomness, nor are we claiming that the only significant developments in algorithmic randomness  are the ones recounted here.  Instead, we aim to provide readers with sufficient context for appreciating the more recent work that is described in the surveys in this volume.  Moreover, we highlight those concepts and results that will be useful for our readers to be aware of as they read the later chapters in this volume.

Before we proceed with the technical material, we briefly comment upon several broader conceptual questions which may occur to the newcomer upon reading this survey:  What is a definition of algorithmic randomness intended to capture?  What is the aim of studying the properties of the various types of randomness?  And why are there so many definitions of randomness to begin with?  It is certainly beyond the scope of this survey to answer these questions in any detail.  Here we note first that more recent motivations for defining randomness and studying the properties of the resulting definitions have become unmoored from the original motivation that led to early definitions of randomness, namely, providing a foundation for probability theory (see, for example, \cite{p14}).  

This original motivation led to the desire for a definition of a random sequence satisfying the standard statistical properties of almost every sequence (such as the strong law of large numbers and the law of the iterated logarithm).  Martin-L\"of's definition was the first to satisfy this constraint.  Moreover, this definition proved to be robust, as it was shown to be equivalent to definitions of randomness with a significantly different informal motivation:  while Martin-L\"of's definition was motivated by the idea that random sequences are statistically typical, later characterizations were given in terms of incompressibility and unpredictability.

With such a robust definition of randomness, one can inquire into just how stable it is:  if we modify a given technical aspect of the definition, is the resulting notion equivalent to Martin-L\"of randomness?  As we will see below, the answer is often negative.  As there are a number of such modifications, we now have a number of nonequivalent definitions of randomness.  Understanding the relationships between these notions of randomness, as well as the properties of the sequences that satisfy them, is certainly an important endeavor.  

One might legitimately express the concern that such work amounts to simply concocting new definitions of randomness and exploring their features.  However, not every new variant of every notion of randomness has proven to be significant.  Typically, attention is given to definitions of randomness that have multiple equivalent formalizations, or which interact nicely with computability-theoretic notions, or which provide insight into some broader phenomenon such as the analysis of almost sure properties that hold in classical mathematics.  Many such developments are outlined in the surveys in this volume.

\subsection{Notation}

Our notation will primarily follow \cite{dhbook} to make it easier to cross-reference these results. The set of natural numbers will be denoted by $\omega$, and we will usually name elements of this set using lowercase Latin letters such as $m$ and $n$. Subsets of $\omega$ will be denoted by capital Latin letters such as $A$ and $B$. Without loss of generality, we may associate an element of $\reals$ (that is, an infinite binary sequence) with the subset of $\omega$ consisting of the places at which the infinite binary sequence is equal to 1. Finite binary strings, or elements of $\strings$, will be denoted by lowercase Greek letters such as $\sigma$ and $\tau$.

We will often wish to discuss the subset of $\reals$ whose elements all begin with the same prefix $\sigma$; we will denote this by $[\sigma]$. We will further extend this to an arbitrary subset $S$ of $\strings$: 
$$[S]=\{ A\in\reals \st \mbox{$\sigma \preceq A$ for some $\sigma\in S$}\}.$$
The first $n$ bits of a binary sequence $X$ of length at least $n$, be it infinite or finite, will be denoted by $X\restr n$, and the length of a finite binary string $\sigma$ will be denoted by $|\sigma|$. The concatenation of two binary strings $\sigma$ and $\tau$ will be denoted by $\sigma\tau$.

\subsection{Computability theory}

This section is intended for researchers in other areas of mathematics who are encountering computability theory for the first time and require an introduction to the underlying concepts; others may safely skip it. While each of Chapter 2 of \cite{dhbook} and Chapter 1 of \cite{niesbook} contains all the fundamental concepts of computability theory that will be required for this volume in more detail, researchers who wish to acquire a deeper understanding of the subject are urged to consult Cooper \cite{cooper04}, Odifreddi \cite{o1,o2}, or Soare \cite{soare}.

Computability theory allows us to think about mathematical functions in an effective context. While there are several ways to formalize the notions we are about to describe, including Turing machines, register machines, the $\lambda$-calculus, and $\mu$-recursive functions, we will not fix such a formalization and will instead encourage the reader to think of the concepts we describe below in terms of the calculations that a computer with potentially unlimited memory is capable of carrying out in a finite but unbounded amount of time.

The most fundamental concept is that of a \emph{partial computable function} $\varphi$, which can be thought of as an idealized computer program that accepts natural numbers as inputs and outputs natural numbers as well. We note that when a computer program is given an input (and thus when a partial computable function is), it may either return an answer at some finite point, or \emph{stage}, or never halt. If a partial computable function halts on every natural number (that is, it is actually a total function), we simply call it a \emph{computable function}. 

We now provide some necessary notation. If $\varphi$ halts on input $n$, we write $\varphi(n)\conv$; if it does not, we write $\varphi(n)\dvrg$. Furthermore, if $\varphi$ halts on input $n$ and gives the output $m$ within $s$ stages, we write $\varphi(n)[s]=m$ to indicate the number of stages as well as the output.

At this point, we make three observations. The first is that there are countably many partial computable functions: each partial computable function is associated with a computer program, and we can note that a computer program is a finite sequence of characters from a finite alphabet and that there are thus countably many such objects. The second is that we can list the partial computable functions in a computable way simply by generating a list of all of the ``grammatically correct" programs and that thus we can speak about, for instance, the $k^{th}$ partial computable function $\varphi_k$. While there are still only countably many (total) computable functions and thus we can list them as well, it can be shown that we cannot list them in a computable way because no computer program is capable of identifying precisely the partial computable functions that halt on every natural number. The third is that there are computable bijections between the natural numbers and the set of finite binary strings and the rational numbers, so we may discuss partial computable functions from or to these sets without loss of generality.

Now we can define special kinds of subsets of $\omega$. A \emph{computably enumerable} (\emph{c.e.}) set is one that is the range of a partial computable function, so we will often write $W_e$ for the set that is the range of the $e^{th}$ partial computable function $\varphi_e$. We may think of $\varphi_e$ as enumerating $W_e$ as follows: we first spend one step trying to compute $\varphi_e(0)$, then two steps trying to calculate each of $\varphi_e(0)$ and $\varphi_e(1)$, then three steps trying to calculate each of $\varphi_e(0)$, $\varphi_e(1)$, and $\varphi_e(2)$, and so on.\footnote{Each of these steps may be said to make up a single stage of the computation mentioned above. These steps will be defined differently based on our formalization of computability theory: they may be the number of states a Turing machine has been in or the number of reduction rules applied in the $\lambda$-calculus, but, at a less formal level, we may think of "spending $n$ steps" as "running the computer program for $n$ seconds."} If the calculation of $\varphi_e(n)$ ever halts, we will eventually discover this through this dovetailing of computations, and when we do, we will enumerate its value into our set $W_e$.  Now we can use this idea of set enumeration to formalize the concept of approximations to a set: for any \ce set $W_e$, we say that the approximation to it at stage $s$ is $W_e[s] = \{n \st  (\exists k\leq s)[\,n=\varphi_e(k)[s]\,]\}$. This gives us a sequence of approximations that converge to our \ce set; in fact, we can write $W_e = \bigcup_{s\in\omega} W_e[s]$. We quickly observe that there are several equivalent definitions of a \ce set; the other one that will be especially useful to us is that of a \ce set as one that is the domain of a partial computable function. 

A set that is itself \ce and has a \ce complement is called a \emph{computable} set. Just as a computable function halts on every input $n$, we can get an answer to ``Is $n$ in $A$?" for a computable set $A$ for every $n$: to see this, we observe that if $A=W_e$ and $\overline{A}=W_i$, then we can determine whether $n\in A$ by enumerating $W_e$ and $W_i$ as described above; $n$ must be in one of them, and we simply note which. These two procedures can again be dovetailed and performed by a single function that will give us the \emph{characteristic function} of $A$:
$$\chi_A(n) = \begin{cases}
1 & n\in A \\
0 & n\not\in A
\end{cases}.$$ 
We can therefore show that the characteristic function of a computable set will be a computable function. Once again, there are countably many \ce sets and countably many computable sets. 

Often, when we discuss randomness, we will talk about a sequence of \emph{uniformly c.e.}\ sets. Instead of simply requiring that each sequence in the set be c.e., we require that there be a single computable function that generates the entire sequence: $\langle A_i\rangle_{i\in\omega}$ is uniformly c.e.\ if there is a computable function $f$ such that $A_i$ is the range of the $f(i)^{th}$ partial computable function. Later, we will generalize this concept to other classes of sets that have some natural indexing: given such a class of sets $\mathcal{C}$, we can say that that we have a sequence of uniformly $\mathcal{C}$ sets $\langle A_i\rangle_{i\in\omega}$  if there is a computable function $f$ such that $f(i)$ gives the index of the $i^{th}$ set in the sequence.

The final topic we must consider in order to understand the concepts in algorithmic randomness we will discuss in this survey is that of oracle computation, or relativization. This requires us to consider \emph{Turing functionals}, usually denoted by capital Greek letters such as $\Phi$, which require not only a natural number $n$ as input but a sequence $X$ that serves as an oracle. These functionals can make use of the standard computational methods of partial computable functions and receive answers to finitely many queries of the sort ``Is $k$ in $X$?" for use in their computation, and they can be indexed as $\Phi_0,\Phi_1,\ldots$  just as the partial computable functions can be indexed as $\varphi_0, \varphi_1,\ldots$. When we use the sequence $X$ as an oracle for the Turing functional $\Phi$, we write $\Phi^X$. Finally, we note that our notation for stages of computations using Turing functionals carries over directly from that for stages of computations using partial computable functions: we write $\Phi_e^X(n)[s]$ just as we would have written $\varphi_e(n)[s]$.

We say that $A$ is \emph{computable from}, or \emph{Turing reducible to} $B$ ($A\leq_T B$) if there is some Turing functional that, given $B$ as an oracle, can compute the characteristic function of $A$. We then use this reducibility to form equivalence classes of sets that we call the \emph{Turing degrees}: $A$ and $B$ have the same Turing degree if $A\leq_T B$ and $B\leq_T A$. This allows us to talk about properties related to a set's computational strength and not its particular members (for instance, we can talk about $\mathbf{0}$, the Turing degree of the computable sets, rather than ``the Turing degree of the empty set"). The Turing degrees will be denoted by boldface lowercase Latin letters such as $\mathbf{d}$.

We will also use relativization to define new sets and Turing degrees. For instance, for each set $A$, we define $A'$ to be $\{n \st \Phi_n^A(n)\conv\}$ and call it the \emph{jump} of $A$. The jump of the empty set, $\emptyset'$, is therefore the set of all natural numbers $n$ such that $\{n \st \Phi_n^\emptyset(n)\conv\}$, or, in other words, the indices of those Turing functions that halt given their own index as input and no additional information. Its Turing degree, $\mathbf{0'}$, is the degree of the famous Halting Problem (see Chapter II.2 of \cite{o1}). Since $A<_T A'$ for any set $A$, we can develop an infinite ascending chain $\mathbf{0} <_T \mathbf{0'}<_T \mathbf{0''} <_T \ldots$ of Turing degrees. We note quickly that in general, the $n^{th}$ jump of $A$ is written as $A^{(n)}$.

Other, stronger reducibilities and their corresponding degree structures have also been found to be useful in the study of randomness. Turing reductions are not required to converge on any input and, when they do converge, the size of the elements of the oracle queried during the computation is not necessarily bounded by any reasonable function. The next type of reducibility is \emph{weak truth-table reducibility}, or $wtt$-reducibility. A $wtt$-reduction $\Phi^A$ is a Turing reduction in which the computation of $\Phi^A(n)$, should it halt, is carried out by querying only the first $f(n)$ bits of $A$ for a given computable function $f$. Finally, the last reducibility we will discuss in this survey is \emph{truth-table reducibility}, or $tt$-reducibility: a $tt$-reduction $\Phi^A$ is one that will converge at every input given any oracle $A$.\footnote{While this was not the original definition of a $tt$-reduction, it is perhaps the most intuitive. The original definition of a $tt$-reduction explains its name and can be found in Chapter III.3 of \cite{o1}.} It can be seen that
$$A\leq_{tt} B \Longrightarrow A\leq_{wtt} B \Longrightarrow A\leq_T B$$
but that none of these implications reverse.

As with Turing reducibility, we can also create equivalence classes of mutually $wtt$- or $tt$-computable sets and study the $wtt$- and $tt$-degrees. We can then ask about the properties of all of these structures---the Turing degrees, for instance, form an upper semilattice---or types of substructures within these structures, such as ideals, which are subsets that are both downward closed and closed under join, or an interval between two degrees (for instance, the interval $[\mathbf{0},\mathbf{0'}]$ in the Turing degrees). We can also discuss relationships between individual degrees; for example, we say that two degrees $\mathbf{c}$ and $\mathbf{d}$ form a \emph{minimal pair} in their degree structure if the only degree they both compute is $\mathbf{0}$.

It is often useful to characterize a subset of $\omega$ in terms of the number of unbounded quantifiers required to define it. A $\Sigma_n$ set $A$ has a membership relation that can be defined from a computable relation $R(x_1,x_2,\ldots,x_n,y)$ using $n$ alternating quantifiers, starting with an existential one: $y\in A$ if and only if 
$$\exists x_1\forall x_2 \exists x_3 \ldots Q_n x_n (R(x_1,x_2,\ldots,x_n,y)).$$
$Q_n$ will be an existential quantifier if $n$ is odd and a universal quantifier if $n$ is even. We observe that we can consider these quantifiers to be strictly alternating. For instance, if we had the membership relation
$$\exists x_1\exists x_2 (R(x_1,x_2,y)),$$
we could use a computable pairing function $p:\omega^2\rightarrow \omega$ and express the same relation as
$$\exists x \exists x_1\leq x \exists x_2\leq x (x=p(x_1,x_2) \wedge R(x_1,x_2,y))$$
instead; note that the second and third existential quantifiers in this formula are bounded and therefore that 
$$\exists x_1\leq x \exists x_2\leq x (x=p(x_1,x_2) \wedge R(x_1,x_2,y))$$
is a computable relation.

\begin{example}
$\emptyset'$ is $\Sigma_1$: $e\in \emptyset'$ if and only if 
$$\exists s (\Phi_e^\emptyset(e)[s]\conv).$$
\end{example}

\begin{example}
The set of all $e$ such that $W_e$ is finite, Fin, is $\Sigma_2$: $e$ belongs to Fin if and only if 
$$\exists m \forall s \forall k (k>m \rightarrow k\not\in W_{e}[s]).$$
\end{example}

A $\Pi_n$ set is defined in a similar way. It will also have $n$ alternating quantifiers, but this time starting with a universal quantifier; we note that the complement of a $\Sigma_n$ set is a $\Pi_n$ set and vice versa. 

\begin{example}
The set of all $e$ such that $\phi_e$ is total, Tot, is $\Pi_2$: $e$ belongs to Tot if and only if 
$$\forall k \exists s \exists m (\phi_e(k)[s]=m).$$
\end{example}

Finally, we have the $\Delta_n$ sets, which we define to be those sets that can be characterized in both a $\Sigma_n$ and a $\Pi_n$ way. Since we often refer to the class of $\Sigma_n$ sets simply as $\Sigma_n$ (and similarly for the classes of $\Pi_n$ and $\Delta_n$ sets), we can write
$$\Delta_n = \Sigma_n \cap \Pi_n.$$

These classes of sets---the $\Sigma_n$, $\Pi_n$, and $\Delta_n$ sets---form the \emph{arithmetic hierarchy}. We will note some fundamental facts relating these classes to the classes of sets we have already discussed (see Chapter IV.1 of \cite{o1}):

\begin{enumerate}
\item $\Sigma_0 = \Pi_0 = \Delta_0 = \Delta_1$ is simply the class of computable sets.
\item A set is $\Sigma_1$ if and only if it is \ce
\item A set is $\Delta_n$ exactly when it is Turing computable from $\emptyset^{(n-1)}$.
\end{enumerate}

Furthermore, this hierarchy is proper: as long as $n>0$, we will always have $\Delta_n\subsetneq \Sigma_n$ and $\Delta_n\subsetneq \Pi_n$.

We may describe subsets of natural numbers in ways other than the arithmetic hierarchy, too. For instance, we have the \emph{high} and \emph{low} sets, which are defined based on the usefulness of the sets in question as oracles:
\begin{itemize}
\item a low set $A\subseteq \omega$ is one such that $A'\equiv_T 0'$, and 
\item a high set $A\subseteq \omega$ is one such that $A'\geq_T 0''$.
\end{itemize}
It is often more useful to define high sets in terms of the domination property discovered by Martin \cite{martin66}: A set is high if and only if it Turing computes a function $f$ that dominates all computable functions (that is, for each computable $g$, we have $f(n)\geq g(n)$ for some sufficiently large $n$). 

Similarly, we may consider highness and lowness in the context of $tt$-reducibility: a set $A$ is \emph{superhigh} if $A'\geq_{tt} 0''$ and \emph{superlow} if $A'\equiv_{tt}0'$.

Another hierarchy of classes of sets that has proven useful is the genericity hierarchy, which we can compare to a hierarchy of randomness notions that we will see later. Instead of classifying sets based directly on the complexity of their definitions, we classify them in terms of the complexity of the sets they are forced to either meet or avoid.

\begin{definition}
Let $S$ be a set of finite binary strings. We say that an infinite binary sequence $A$ \emph{meets} $S$ if there is some $\sigma\in S$ that is an initial segment of $A$. Furthermore, we say that $A$ \emph{avoids} $S$ if there is some initial segment of $A$ that is not extended by any element of $S$.
\end{definition}

This gives us the framework necessary to define generic sets once we have also defined a dense set of strings: a set of strings $S$ is \emph{dense} if for every $\tau\in\strings$, there is a $\sigma\in S$ that extends it.

\begin{definition}
A sequence $A$ is \emph{$n$-generic} if it either meets or avoids every $\Sigma_n$ set and \emph{weakly $n$-generic} if it meets every dense $\Sigma_n$ set.
\end{definition}

These classes of sequences once again form a proper hierarchy: every $n$ generic is weakly $n$-generic, and every weakly $(n+1)$-generic is $n$-generic.

Other classes of sets whose definitions are less closely tied to the arithmetical hierarchy have also been shown to be useful. For instance, we will make use of the sets of hyperimmune degree, which are defined, once again, using a domination property \cite{mm68}: a set $A$ has \emph{hyperimmune degree} if it computes a function that is not dominated by \emph{any} computable function. The sets that do not have hyperimmune degree are said to be of \emph{hyperimmune-free degree}; it is worth noting that there are continuum many such sets and that all of them (except the computable sets) are Turing incomparable to $\emptyset'$. 

\subsection{Core mathematical concepts}\label{ssec:core}

Several concepts from classical mathematics will prove useful; we summarize them here. First we will recall some fundamental facts about the Cantor space, $\reals$, as a topological space and as a probability space.

In the Cantor space, our basic open sets have the form $[\sigma]$ for $\sigma\in\strings$: as previously mentioned, $[\sigma]$ is the set of elements of $\reals$ that extend $\sigma$. In fact, these sets are all clopen, and the clopen sets are the finite unions of these $[\sigma]$'s. 
Now that we have done this, we can describe the complexity of the open sets that we generate in this way in terms of their generating sets. For instance, we can say that $[S]$ is effectively open if $S$ is c.e., and more generally, we can define the effective Borel hierarchy as we defined the arithmetic hierarchy in the previous subsection: $[S]$ is $\Sigma^0_1$ if it is the union of a \ce sequence of basic open sets, $\Pi^0_1$ if it is the complement of a $\Sigma^0_1$ set, $\Sigma^0_n$ for $n>1$ if it is the union of a  sequence of uniformly $\Pi^0_{n-1}$ sets (that is, a computable sequence of $\Pi^0_{n-1}$ classes), and so on. It is worth noting at this point that it is customary to refer to a subset of the Cantor space, especially one defined using this hierarchy, as a class. 

We can also establish the Lebesgue measure on the Cantor space: the measure of a basic open set $[\sigma]$ is $\mu([\sigma])=2^{-|\sigma|}$, and the measure of any other measurable set is determined in the standard way.

We will often identify the Cantor space with the unit interval (0,1) since these spaces are measure-theoretically isomorphic. Here, we make use of the interval topology on $\R$, and our basic open sets are intervals $(a,b)$. We will establish the Lebesgue measure in this context as well, denoted throughout by $\mu$ once again \cite{lebesgue}. This is the ``standard" measure on $\R$, and the Lebesgue measure of such an interval is $\mu((a,b))=b-a$ for finite $a$ and $b$. 

We can also describe elements of $\R$ using concepts from classical computability theory. In general, we identify a real $\alpha$ in the unit interval with the element $A$ of the Cantor space such that $\alpha=0.A$.\footnote{While some reals may have two representations, this does not matter: such a real will be rational and therefore the corresponding possibilities for $A$ are computable and will both have the same computational strength.}  This allows us to say that such a real is computable if the corresponding $A\in\reals$ is; it is equivalent to say that a real $\alpha$ is computable if there is a computable sequence of rationals $\langle q_i\rangle$ converging to it such that $|q_n-\alpha|<2^{-n}$ for every $n$ \cite{turing36,turing37}. 

Of course, we would like to extend this to computable enumerability as well: just as a \ce set is one which we build up from $\emptyset$ by enumerating elements into it, a \emph{left-c.e.\ real} $\alpha$ is one that is effectively approximable from below; that is, there is a computable, increasing sequence of rationals that converges to $\alpha$. Correspondingly, a \emph{right-c.e.\ real} is one that is effectively approximable from above. Equivalently, we could define these in terms of Dedekind cuts: a real $\alpha$ is left-\ce if and only if its left cut $\{q\in\Q \st q<\alpha\}$ is a c.e.\ set and right-\ce if its right cut (defined similarly) is a \ce set.

We can extend the notions of computability and computable enumerability once more, this time to functions from $\omega$ or another computable set to $\R$. To do so, we need to be able to talk about a sequence of reals that is uniformly computable or left-c.e. These definitions are built directly from those of uniformly computable and uniformly \ce sets:

\begin{definition}
A uniformly computable (left-c.e.) sequence of reals is a sequence $\langle r_i\rangle_{i\in\omega}$ such that there is a computable function $f:\omega^2\rightarrow \Q$ such that for a given $i$, $\langle f(i,n)\rangle_{n\in\omega}$ is a computable (left-c.e.) approximation for $r_i$.
\end{definition}

\begin{definition}
A function from a computable set to $\R$ is computable if its values are uniformly computable reals, and it is computably enumerable if its values are uniformly left-c.e.\ reals.
\end{definition}

Now we turn our attention to the general mathematical ideas we will need to study randomness properly and place them in the context of computability theory. The first concept, that of a martingale, will be useful when we discuss the predictability framework for randomness. In general, a martingale is a certain type of stochastic process, but here we need only think of it as a type of betting strategy on finite binary strings.

\begin{definition}\cite{levy37}
A function $d:\strings\rightarrow \R^{\geq 0}$ is a \emph{martingale} if it obeys the fairness condition 
$$d(\sigma)=\frac{d(\sigma0)+d(\sigma1)}{2}$$
for all $\sigma\in\strings$. We say that a martingale $d$ \emph{succeeds} on $A\in\reals$ if 
$$\limsup_n d(A\restr n) = \infty,$$
and the \emph{success set} of $d$, which we will write as $S[d]$, is the set of all sequences upon which $d$ succeeds.
\end{definition}

We can think of $d(\sigma)$ as expressing the amount of capital that we have after betting on the initial string $\sigma$ using the strategy inherent in $d$ (so $d(\langle\rangle)$ is our  capital before any bets are placed) and $S[d]$ as the set of sequences that we can make arbitrarily much money betting on if the payout is determined by $d$.  As shown by Ville,  $P\subseteq2^\omega$ has Lebesgue measure zero if and only if there is some martingale $d$ such that $P\subseteq S[d]$ \cite{ville}.

A computable or \ce martingale is simply a martingale that is, respectively, a computable or \ce function.

To define the last of the core mathematical concepts in this section, Hausdorff dimension, we must consider a variation of Lebesgue measure on the Cantor space. 

\begin{definition} \cite{hausdorff18}
Let $0\leq s\leq 1$. The \emph{$s$-measure} of a basic open set $[\sigma]$ is $\mu_s([\sigma]) = 2^{-s|\sigma|}$.
\end{definition}

This will allow us to define dimensions of subsets of Cantor space.

\begin{definition}
An \emph{$n$-cover} of a subset $S$ of Cantor space is a set of strings $C\subseteq 2^{\geq n}$ such that $S\subseteq [C]$.  We define
$$H^s_n(S) = \inf \left\lbrace \sum_{\sigma\in C} \mu_s([\sigma]) \st \mbox{$C$ is an $n$-cover of $S$}\right\rbrace.$$

The \emph{$s$-dimensional outer Hausdorff measure} of $S$ is 
$$H^s(S) = \lim_n H^s_n(S),$$
and the \emph{Hausdorff dimension} of $S$ is
$$\dim_H(S) = \inf\{ s\st H^s(S)=0\}.$$
\end{definition}
\noindent  We note that these definitions are specialized here to Cantor space. For details about these notions in more general settings, see, for instance, the monographs \cite{rogers70} or  \cite{falconer14}.

We further note that the effective version of Hausdorff dimension that we will discuss later, originally due to Lutz (see, e.g. \cite{lutz03}), is instead based on a type of betting strategy called an $s$-gale rather than $n$-covers.  In fact, Lutz provided an equivalent characterization of classical dimension in terms of $s$-gales and then gave an effectivization of this alternative notion.

\section{Early developments}\label{sec:early}

During the 1960s and early 1970s, the foundation of much of the current work in algorithmic randomness was laid.  Earlier work by von Mises, Ville, Church, and Wald in the first half of the twentieth century highlighted the problem of defining the notion of an individual random sequence, but no satisfactory definition of a random sequence was provided during this time.  However, several promising definitions of randomness emerged in the work of (i) Kolmogorov in the mid-1960s, (ii) Martin-L\"of in the late 1960s, and (iii) Schnorr and Levin in the early 1970s.

We will not rehearse the broader philosophical concerns that motivated the search for a definition of an individual random sequence (for such an account, see \cite{p14}).  For our purposes, it suffices to note the key desiderata for a definition of randomness that came into focus during the first half of the twentieth century:  (1) a random sequence is one that should not be contained in any null sets that are ``nicely" definable in some way (clarified in part by the work of Wald and Ville on von Mises' original definition of randomness) and (2) one should formalize these ``nicely" definable null sets in terms of effectivity (as suggested by Church's introduction of concepts of computability theory to the task of defining randomness).

\subsection{Randomness via initial segment complexity}\label{subsec-isc}

The first significant breakthroughs in the task of offering such a definition of randomness came in the work of Kolmogorov and, independently, Solomonoff \cite{solo64a,solo64b},  who provided different accounts of the initial segment complexity of sequences.  We will focus here on Kolmogorov's contribution.

Kolmogorov did not set out to provide a definition of random sequences in terms of some class of effective null sets.  Rather, his aim was to provide a notion of randomness for finite strings (again, for the motivation behind this aim, see \cite{p14}), which was, in turn, defined in terms of Kolmogorov complexity.  Such a definition is found in Kolmogorov's  1965 paper ``Three approaches to the quantitative definition of information" (see \cite{kol1} for the original English translation of the article).\footnote{Here we do not exactly follow the details of Kolmogorov's presentation; he initially defines conditional Kolmogorov complexity in his paper.}  For a fixed partial computable function $M: \strings\rightarrow\strings$, often called a \emph{machine}, and some $\tau\in\strings$, we can consider all strings $\sigma\in\strings$ such that $M(\sigma){\downarrow}=\tau$; if we consider each such string $\sigma$ to contain information about how to produce the string $\tau$ (via the function $M$), then the shortest such string $\sigma$ provides the minimal amount of information necessary for producing $\tau$.  This gives us a measure of the complexity of $\tau$:

\begin{definition}
The \emph{plain Kolmogorov complexity} of a string $\tau$ relative to the machine $M$ is
\[
C_M(\tau)=\min\{|\sigma|\st M(\sigma){\downarrow}=\tau\},
\]
assuming there is some $\tau\in\strings$ such that $M(\sigma){\downarrow}=\tau$; if no such $\tau$ exists, we set $C_M(\tau)=\infty$.  
\end{definition}

The dependence of the above definition on the function $M$ is undesirable, for clearly, different choices of $M$ produce different values $C_M(\tau)$ for each fixed $\tau\in\strings$.  Kolmogorov addressed this problem by defining his complexity measure in terms of a \emph{universal} partial computable function:  If $\langle M_e\rangle_{e\in\omega}$ is a fixed computable enumeration of partial computable functions from $\strings$ to $\strings$, then we can define a universal partial computable $U:\strings\rightarrow\strings$  by setting $U(1^e0\sigma)=M_e(\sigma)$ for each $e\in\omega$ and $\sigma\in\strings$ (assuming that $M_e(\sigma){\downarrow}$; otherwise, we set $U(1^e0\sigma){\uparrow}$).  Then we define $C(\sigma)$ to be $C_U(\sigma)$ for each $\sigma\in\strings$.

The value in defining complexity in terms of a universal partial computable function is seen in the following result, which is often referred to as the invariance theorem:

\begin{thm}[Kolmogorov, \cite{kol1}]
For every partial computable $M:\strings\rightarrow\strings$, there is some $c\in\omega$ such that for every $\sigma\in\strings$,
\[
C(\sigma)\leq C_M(\sigma)+c.
\]
\end{thm}

Note that the invariance theorem implies that the complexity values determined by two different choices of universal partial computable functions $U$ and $U'$ yield complexity measures that only differ by a finite fixed constant $c$:  $|C_U(\sigma)-C_{U'}(\sigma)|\leq c$ for every $\sigma\in\strings$. 

From Kolmogorov's notion of complexity, how do we define the randomness of binary strings?  
Let us first consider two informal examples, a string $\sigma_1$ consisting of 50,000 1's and a string $\sigma_2$ of length 50,000 obtained by the tosses of a fair coin.  The shortest program needed to generate $\sigma_1$ has length considerably shorter than 50,000, as such a program only needs to specify that the symbol `1' is to be repeated 50,000 times.  However, any program that generates $\sigma_2$ must, with high probability, have most, if not all, of the entire string $\sigma_2$ hardwired into the program, as most strings obtained by tossing a fair coin contain very few regularities and thus cannot be compressed.  Following this example, the idea behind Kolmogorov's definition of randomness is to identify randomness with incompressibility.  We now turn to the formal details.

Kolmogorov did not explicitly define randomness in his 1965 paper, but we find the following definition in the 1969 paper \cite{kol69}.   Consider those strings $\sigma$ such that $C(\sigma)\geq |\sigma|$.  As there is no program to generate such a string $\sigma$ that is shorter than $\sigma$, the most efficient way to produce such a string via $U$ is simply to give it as input and copy it directly to output.

More generally, for a fixed $c\in\omega$, we can consider the set of all strings that cannot be compressed by more than $c$ bits.  Let us say of a string $\sigma\in\strings$ that it is $c$-\emph{incompressible} if $C(\sigma)\geq |\sigma|-c$.  For $n\geq c$, since the number of strings of length strictly less than $n-c$ is equal to $1+2^1+2^2+\dotsc+2^{n-c-1}=2^{n-c}-1$, there are at least $2^n(1-2^{-c})-1$ $c$-incompressible strings of length $n$, thereby yielding a plethora of random strings (for sufficiently large $n$).

\subsection{\ml randomness}\label{subsec-ml}
It is natural to try to extend this definition of randomness for finite binary strings to infinite binary sequences, a task that Martin-L\"of sought to accomplish.  Indeed, we can define a sequence $A\in\reals$ to be $c$-incompressible if there is some $c\in\omega$ such that $C(A\restr n)\geq n-c$ for every $n\in\omega$.  However, as shown by Martin-L\"of (see \cite[Theorem 3.1.4]{dhbook}), no sequence is $c$-incompressible for every $c\in\omega$.  In fact, he shows that for every $A\in\reals$, there are infinitely many $n$ such that $C(A\restr n)\leq n-\log(n)$ (where $\log(n)$ is the binary logarithm).  

Instead of seeking to define random infinite sequences by modifying the underlying notion of Kolmogorov complexity, Martin-L\"of in \cite{ml66} took an alternative approach, defining randomness in terms of certain effective statistical tests.  Martin-L\"of aimed to formalize the notion of statistical tests used in hypothesis testing, where each such test has a critical region such that, given a sequence of observations contained in this region, we reject the null hypothesis at a certain level of significance.  In Martin-L\"of's formalization,  such a test is given by a computable sequence of effectively open sets $\langle U_i\rangle_{i\in\omega}$, each of which corresponds to a critical region given by a certain level of significance.   Moreover, if a sequence is random, then it should not be contained in the critical regions at every level of significance; eventually, we should find a critical region at some level of significance that the random sequence avoids.  Thus we have the following definition:

\begin{definition}
A \emph{\ml test} is a sequence of uniformly $\Sigma^0_1$ classes $\langle U_i\rangle_{i\in\omega}$ such that $\mu(U_i)\leq 2^{-i}$ for every $i\in\omega$. A sequence $A\in\reals$ is said to \emph{pass} a Martin-L\"of test $\langle U_i\rangle_{i\in\omega}$ if $A\notin\bigcap_{i\in\omega} U_i$, and a sequence is \emph{Martin-L\"of random} if it passes every Martin-L\"of test.
\end{definition}

Martin-L\"of also proved in \cite{ml66} that there is a \emph{universal} Martin-L\"of test: a test $\langle U_i\rangle_{i\in\omega}$ such that for all Martin-L\"of tests $\langle V_i\rangle_{i\in\omega}$, $\bigcap_{i\in\omega} V_i\subseteq \bigcap_{i\in\omega} U_i$. This provides an easy argument that the class of Martin-L\"of random sequences has measure 1: since the sequences that are not Martin-L\"of random are precisely those that are not contained in $\bigcap_{i\in\omega} U_i$ for a universal test $\langle U_i\rangle_{i\in\omega}$, the class of Martin-L\"of random sequences must be conull. Moreover, Martin-L\"of showed that common laws of probability such as the strong law of large numbers and the law of the iterated logarithm are satisfied by all Martin-L\"of random sequences.

Martin-L\"of also investigated the extent to which one can characterize Martin-L\"of randomness in terms of initial segment complexity, showing that for every sequence $A$ such that, for some $c\in\omega$,
$C(A\restr n)\geq n-c$ holds for infinitely many $n\in\omega$, $A$ is Martin-L\"of random (see \cite{ml71}).  He was, however, unable to establish the converse, which in fact does not hold, as shown much later independently by J.\ Miller \cite{miller04} and Nies, Stephan, and Terwijn \cite{nst05} 
(see Section \ref{subsec-lowness} below for more details).

The problem of providing an initial segment complexity characterization of randomness for infinite sequences was solved independently by Levin and Schnorr, who worked with an alternative notion of Kolmogorov complexity, namely, prefix-free Kolmogorov complexity.  Recall that a set $S\subseteq\strings$ is \emph{prefix-free} if for every $\sigma\in S$, if $\tau\succ\sigma$, then $\tau\notin S$.  We can extend the notion of being prefix-free to a machine $M:\strings\rightarrow\strings$ by stipulating that $M$ is prefix-free if the domain of $M$ is a prefix-free subset of $\strings$.  Both Levin \cite{levin73} and Chaitin \cite{c75} independently defined prefix-free Kolmogorov complexity in terms of a prefix-free universal machine.  Note that one can effectively enumerate the collection of all prefix-free machines $\langle M_e\rangle_{e\in\omega}$ and hence we can define $U$ as we did in Subsection \ref{subsec-isc} above.

\begin{definition}
The \emph{prefix-free Kolmogorov complexity} of a string $\tau$ is given by
\[
K(\tau)=\min\{|\sigma|\st U(\sigma){\downarrow}=\tau\},
\]
where $U$ is a universal, prefix-free machine.
\end{definition}

For a  prefix-free machine  $M$, we write $K_M$ to stand for prefix-free complexity relative to the machine $M$. One can readily verify that results such as the invariance theorem and the existence of incompressible strings still hold for prefix-free complexity.

With this modified notion of Kolmogorov complexity, Levin and Schnorr independently proved the following:\footnote{Schnorr is given credit for this result in Chaitin's \cite{c75}.}

\begin{thm}[Levin \cite{levin73}, Schnorr \cite{c75}]
$A\in\reals$ is Martin-L\"of random if and only if there is some $c\in\omega$ such that for all $n\in\omega$,
\[
K(A\restr n)\geq n-c.
\]
\end{thm}

As part of the proof of the Levin-Schnorr theorem, it is standard to prove that for every $A$ that is not Martin-L\"of random, there is some prefix-free machine $M$ such that for every $c\in\omega$, there is an $n\in\omega$ such that $K_M(A\restr n)<n-c$ (where $K_M$ is prefix-free complexity relative to the machine $M$).  More specifically, given a Martin-L\"of test that $A$ fails to pass, we explicitly construct such a machine $M$.  This is done via what was earlier referred to as the Kraft-Chaitin theorem but is currently called the KC-theorem (as in \cite{dhbook}) or the machine existence theorem (as in \cite{niesbook}).  The theorem is as follows; see the footnote in Section 3.6 of \cite{dhbook} for an explanation of its origins.

\begin{thm}
Given a computable enumeration of pairs $\langle(n_i,\sigma_i)\rangle\subseteq \omega\times\strings$ such that $\sum_{i\in\omega}2^{-n_i}\leq 1$, there is a prefix-free machine $M$ and a sequence of strings $\langle \tau_i\rangle$ such that for every $i$, $M(\tau_i)=\sigma_i$ and $|\tau_i|=n_i$.
\end{thm}

To see the details of how this result is used in the proof of the Levin-Schnorr theorem, see \cite[Theorem 3.2.9]{niesbook}.

The Levin-Schnorr theorem thus establishes the equivalence of a characterization of randomness in terms of statistical properties with a characterization of randomness in terms of incompressibility.  In fact, the proof of the Levin-Schnorr theorem shows that (i) if a sequence is compressible, then there is a statistical test such that its critical regions contain the sequence at all levels of significance (each level corresponding to how compressible the sequence is), and (ii) if a sequence is statistically atypical, contained in all critical regions at every level of significance of some statistical test, then we can define a machine to compress this sequence (and in fact, every sequence contained in all such critical regions).
 
\subsection{Schnorr's contributions}

Shortly after the publication of Martin-L\"of's definition of randomness, Schnorr developed an alternative approach to algorithmic randomness using martingales as originally defined by Ville.  This approach constitutes the third major framework for algorithmic randomness: the unpredictability framework.  Not only did Schnorr provide a characterization of Martin-L\"of randomness in terms of unpredictability, but he also introduced two alternative definitions of randomness, computable randomness and Schnorr randomness in \cite{schnorr71} and \cite{schnorr}; these definitions make use of computable and \ce martingales as described in Section \ref{ssec:core}.

\begin{thm}[Schnorr, \cite{schnorr71,schnorr}]
A sequence is \ml random if and only if no \ce martingale succeeds on it.
\end{thm}

In other words, a sequence is \ml random exactly when no \ce martingale can predict it sufficiently well to make an arbitrarily large amount of money on it. We quickly note that the standard proof of this theorem involves interpreting a \ml test as a \ce martingale and that therefore, since there is a universal \ml test, there is also a universal \ce martingale.

As noted above, in the same works, we find two new randomness notions, namely, computable randomness and the notion now known as Schnorr randomness. The definition of Schnorr randomness requires the definition of an \emph{order function}, a nondecreasing and unbounded computable function from $\omega$ to $\omega$.\footnote{It is usual at the time of writing to require only that an order function be increasing and unbounded. For instance, this is the definition used in Downey and Hirschfeldt \cite{dhbook}. However, when Schnorr originally defined order functions, he required that they be computable as well, and since we will not consider noncomputable order functions in this survey, we prefer to use his original definition for simplicity's sake.} 

\begin{definition}
A sequence is computably random if no computable martingale succeeds on it.
\end{definition}

\begin{definition}
A sequence $A$ is Schnorr random if, for every computable martingale $d$ and order function $p$, 
$$\limsup_n \frac{d(A\restr n)}{p(n)} < \infty.$$
\end{definition}

In other words, a sequence is computably random if no computable martingale can predict it well enough to make an arbitrarily large amount of money betting on it, and a sequence is Schnorr random if, while a computable martingale may be able to predict it well enough to make arbitrarily much betting on it, we can bound the rate at which this occurs.  

In \cite{schnorr}, Schnorr also showed that Schnorr randomness could be characterized using a more effective form of \ml tests: rather than require the measure of each test component to be effectively approximable from below, one should require it to be approximable from both below and above.

\begin{definition}
A Schnorr test is a \ml test $\langle V_i\rangle_{i\in\omega}$ such that the measure of the $i^{th}$ test component, $\mu(V_i)$, is exactly $2^{-i}$. 
\end{definition}

We note here that it is actually only necessary that the measures of the test components be uniformly computable; however, it is customary to use $2^{-i}$ as the measure of the $i^{th}$ component. 

The primary difference between Martin-L\"of randomness and the two alternative definitions that Schnorr proposed is that the latter two definitions are more constructive than the former.  First, Schnorr objected to the notion of a c.e.\ martingale, questioning why one should only require that the values of a martingale be computably approximable from below but not from above.  Thus, in the move from c.e.\ martingales to computable martingales, we arrive at the definition of computable randomness.  Second, Schnorr rejected the standard notion of martingale success, pointing out the possibility of winnings of a gambler increasing so slowly that this success goes undetected.  A better approach, on Schnorr's view, is to require a computable lower bound on the success of a martingale, so that, at least in principle, the gambler can recognize that her winnings are growing without bound.  In the move from the standard notion of success to his more constructive alternative, we then arrive at the definition of Schnorr randomness.

While Schnorr randomness and computable randomness were the first variations on \ml randomness to be considered, many others have been presented since then. While several will be discussed in later sections, two general approaches will not appear in this survey. One is higher randomness, which allows randomness notions defined in terms of effective descriptive set theory and is discussed in Monin's survey. Another is resource-bounded randomness, which allows randomness notions whose definitions involve time and space bounds and is discussed in Stull's survey.

\section{Intermittent work:  The late twentieth century}\label{sec:middle}

After the initial flurry of results established primarily by Martin-L\"of, Schnorr, and Levin in the late 1960s and early 1970s, intermittent work on algorithmic randomness was carried out in the mid-1970s, the 1980s, and the early 1990s.  We highlight some of these developments, starting with the work of Demuth and Ku\v cera.

\subsection{The contributions of Demuth and Ku\v cera}\label{subsec:demkuc}
Working in the tradition of the Markov school of constructive analysis, the Czech mathematician Osvald Demuth made a number of important contributions to the study of algorithmic randomness.  As laid out in the survey \cite{dempath}, not only did Demuth independently discover the notions of Martin-L\"of randomness and computable randomness, but he also developed several notions of randomness that have proven to be fruitful in the decades since then. For a survey on Demuth's use of notions of randomness in constructive analysis, particularly the notions now called Demuth randomness and weak Demuth randomness, see \cite{dempath}.  For more on Demuth's contributions to computable analysis, see both the Rute survey and the Porter survey in this volume.

Ku\v cera, who collaborated with Demuth on a number of projects in computable analysis and effective notions of genericity, also proved a number of key results in algorithmic randomness in the mid-1980s.  We focus primarily on the paper \cite{k85}, which contains several results that are now considered to be classical theorems in algorithmic randomness.  

\begin{thm}[Ku\v cera \cite{k85}]
For every $\Pi^0_1$ class $\mathcal{P}$ with $\mu(\mathcal{P})>0$ and every Martin-L\"of random $A\in\reals$, there is some tail $B$ of $A$, that is, a $B$ such that $A=\sigma B$ for some $\sigma\in\strings$, such that $B\in\mathcal{P}$.  In particular, the collection of Turing degrees of members of $\mathcal{P}$ includes every degree containing a Martin-L\"of random sequence.
\end{thm}

As noted in \cite{bdms10}, another way to formulate the first part of this theorem is that for any $\Pi^0_1$ $\mathcal{P}$ class of positive measure and any Martin-L\"of random sequence $X$, we can apply the shift operator (each application of which drops the initial bit of a sequence) to $X$ a finite number of times to eventually obtain an element of $\mathcal{P}$.  Seen in this light, this result relates ergodic theory and algorithmic randomness, a topic that has been fruitfully explored, as laid out in Towsner's survey in this volume.

Next, recall that a sequence $A$ has \emph{diagonally noncomputable degree} (or \emph{DNC degree}) if there is some $f\leq_T A$ such that $f(i)\neq \phi_i(i)$ for every $i\in\omega$.  Sequences of DNC degree are closely related to completions of Peano arithmetic, which are studied in the context of the G\"odel incompleteness phenomenon.  In this same paper \cite{k85}, Ku\v cera proved the following (in a slightly different form):

\begin{thm}
Every Martin-L\"of random sequence has DNC degree.
\end{thm}

Lastly, we also find the following in \cite{k85}:

\begin{thm}
For every Turing degree $\mathbf{a} \geq\mathbf{0'}$, there is some Martin-L\"of random $A\in\mathbf{a}$.
\end{thm}

Implicit in the proof of this last result is the theorem that for every $B\in\reals$, there is some Martin-L\"of random $A\in\reals$ such that $B\leq_T A$.  This result is now referred to as the Ku\v cera-G\'acs theorem, as it was independently established by G\'acs \cite{g86}, who actually showed that there is a computable bound on the use of the computation of $B$ from $A$, in other words, that $B\leq_{\mathrm{wtt}}A$.  

Ku\v cera's result is counterintuitive:  How can every sequence that computes the halting problem be Turing equivalent to a \ml random sequence?  It turns out that this phenomenon is the exception rather than the rule, for as Sacks' Theorem states \cite{sacksdou}, for every $A\in\reals$, there are measure zero many sequences $X\in\reals$  such that $A\leq_TX$.  Thus, only measure zero many Martin-L\"of random sequences are computationally strong enough to compute the halting problem.  Moreover, as we will see in Section \ref{subsec:randseqTdegs} below, this phenomenon is incompatible with stronger notions of randomness.
 
\subsection{The contributions of Kurtz, Kautz, and van Lambalgen}\label{subsec-kkvl}
Now we turn to the work of Kurtz, Kautz and van Lambalgen. While Ku\v{c}era focused on situating the \ml random sequences within the Turing degrees, Kurtz created true hierarchies of randomness notions for the first time in his 1981 dissertation \cite{kurtz}: the $n$-random and weakly $n$-random sequences. 

\begin{definition}
$A\in\reals$ is \emph{$n$-random} if for every sequence $\langle S_i\rangle_{i\in\omega}$ of uniformly $\Sigma^0_n$ subsets of $\reals$ such that $\sum_{i\in\omega} \mu(S_i)$ is finite, $A$ belongs to only finitely many $S_i$'s. $A$ is \emph{weakly $n$-random} if it is not contained in any $\Pi^0_n$ class of Lebesgue measure zero, or equivalently, if it belongs to every conull $\Sigma^0_n$ subset of $\reals$.
\end{definition}

The relationships between these hierarchies that we would expect to hold based on their names actually do: every $n$-random set is weakly $n$-random, every weakly $(n+1)$-random set is $n$-random, and these inclusions are strict for all $n\geq 1$.\footnote{Kurtz himself only proved that there are $n$-random degrees that are not weakly $(n+1)$-random; the other separation was proven later by Kautz \cite{kautz}.} We can also note that a result of Solovay shows us that the 1-random sets are precisely the \ml random sets \cite{solovay75}, giving us another type of test characterization for \ml randomness.  Furthermore, weak 1-randomness has been extensively studied in its own right and is now known as \emph{Kurtz randomness}. For instance, Kurtz showed that every hyperimmune degree is in fact Kurtz random. This notion of randomness is therefore very weak, and it can reasonably be claimed that it is not a ``proper" randomness notion.\footnote{We note that randomness notions that result in large classes of random sequences are often called \emph{weak}, even though the tests used to describe such notions are necessarily more restricted: it is easier to pass every test in a smaller class of tests. Similarly, a randomness notion describing a relatively small class of random sequences is called \emph{strong}.}
 
Ten years later, Kautz would continue Kurtz's investigations of $n$-randomness and weak $n$-randomness as well as his own measure-theoretic pursuits in his dissertation \cite{kautz}. In particular, he showed that $n$-randomness can be characterized using more traditional methods such as tests of the same form as \ml tests, Kolmogorov complexity, and martingales. We define a \emph{$\Sigma^0_n$-test} to be a sequence $\langle V_i\rangle_{i\in\omega}$ of uniformly $\Sigma^0_n$ classes such that $\mu(V_i)\leq 2^{-i}$ for all $i$. We also consider randomness relative to an oracle for the first time in this survey. When we say that $A$ is random relative to $B$ (for a given randomness notion), we mean in general that access to $B$ does not allow us to derandomize $A$. For $n$-randomness, this is straightforward: we simply use $\emptyset^{(n-1)}$ as an oracle for the universal machine, the universal martingale, or the components of the universal test that we are using to determine the randomness of $A$.     However, for other randomness notions, it is far less straightforward; see Franklin's survey in this volume for a discussion of the issues that arise.

We will write $K^A(\sigma)$ as the prefix-free Kolmogorov complexity of $\sigma$ relative to $A$.
 
 \begin{thm}[Kautz, \cite{kautz}]
For $A\in\reals$, the following are equivalent:
\begin{enumerate}
\item $A$ is $n$-random.
\item For every $\Sigma^0_n$ test $\langle V_i\rangle$, $A\not\in \bigcap_i V_i$.
\item $A$ is 1-random relative to $\emptyset^{(n-1)}$.
\item For some $c\nat$, $K^{0^{(n-1)}}(A\restr m) \geq m-c$ for all $m$.
\item No $\Sigma^0_n$ martingale succeeds on $A$.
\end{enumerate}
\end{thm}  
 
Kurtz also studied the Turing jumps of $n$-randoms; for instance, he shows that the class $\{A\st A^{(n-1)} \geq_T \emptyset^{(n)}\}$ is null and does not contain any $(n+1)$-random sets. 

Finally, Kautz considered randomness with respect to computable measures in general and not just the Lebesgue measure, allowing him to define $n$-randomness with respect to an arbitrary computable measure $\nu$. He further examined the extent to which 1-randomness is preserved under effective transformations, and how to translate between randomness with respect to various computable measures. Porter continues this discussion of randomness with respect to different measures, both computable and noncomputable, in his survey.
  
Kautz also demonstrated that any computable subset of an $n$-random sequence is itself $n$-random and that, in particular, a computable subset of an $n$-random sequence cannot be computed from ``the rest" of the sequence and is in fact $n$-random relative to it \cite{kautz}. In contrast, he also showed that for every $n\geq 1$, there is a weakly $n$-random set $A\oplus B$ such that $A$ is not weakly $n$-random relative to $B$. This work is parallel to that of van Lambalgen who, in \cite{vl90}, precisely characterized the circumstances under which the join of two \ml random sequences would itself be random:

\begin{thm}[van Lambalgen, \cite{vl90}]
For $A,B\in\reals$, the following are equivalent:
\begin{enumerate}
\item $A\oplus B$ is $n$-random.
\item $A$ is $n$-random and $B$ is $n$-random relative to $A$.
\end{enumerate}
\end{thm}

As an immediate consequence, we get that a sequence $A$ is $n$-random relative to an $n$-random sequence $B$ if and only if $B$ is $n$-random relative to $A$.  This symmetry of relative randomness has proven to be an extremely useful tool in the study of Martin-L\"of randomness.  A more detailed discussion of this theorem, specifically in the context of alternative notions of randomness, can be found in Franklin's survey in this volume.

\section{Rapid growth at the turn of the century}\label{sec:naughties}

During the early 2000s, interest in algorithmic randomness grew considerably in the computability theory community.  Here we highlight some of the more significant results. 

\subsection{The Turing degrees of random sequences}\label{subsec:randseqTdegs}

As discussed in Section \ref{sec:early}, Schnorr introduced two alternatives to Martin-L\"of randomness, namely Schnorr randomness and computable randomness.  Separations of these notions were initially established by Schnorr \cite{schnorr71} and Wang \cite{wang, wang99}.  The results were extended by Nies, Stephan, and Terwijn \cite{nst05}, who proved that these notions could be separated in precisely the high degrees.

\begin{thm} [Nies, Stephan, and Terwijn \cite{nst05}] The following are equivalent for $A\in\reals$:
\begin{itemize}
\item[(i)] $A$ has high Turing degree (that is, $A'\geq_T\emptyset''$).
\item[(ii)] There is some $B\equiv_TA$ that is computably random but not Martin-L\"of random.
\item[(iii)] There is some $C\equiv_TA$ that is Schnorr random but not computably random.
\end{itemize}
\end{thm}

Nies, Stephan, and Terwijn also proved that this theorem holds when we consider only left-c.e.\ reals, so we can separate these notions even in that more limited context. 

Another degree-theoretic result established in this period by Stephan characterizes the Martin-L\"of random sequences that are Turing complete.

\begin{thm}[Stephan \cite{stephan02}]\label{stephan-thm}
For every Martin-L\"of random sequence $A\in\reals$, $A\geq_T\emptyset'$ if and only if $A$ has PA degree, that is, $A$ computes a consistent completion of Peano arithmetic.
\end{thm}

Recall from the previous section that Ku\v cera proved that every Turing degree above $\mathbf{0}'$ contains a Martin-L\"of random sequence.  From Stephan's results, these are precisely the PA degrees that contain a Martin-L\"of sequence. Later, Franklin and Ng would characterize these using a new randomness notion, difference randomness \cite{fng}.

Lastly, a key result on the Turing degrees of Martin-L\"of random sequences relative to various oracles, known as the ${XY\!Z}$ theorem, is the following:

\begin{thm}[J.\ Miller and Yu \cite{milleryu08}]
For a Martin-L\"of random sequence $X$, if $X\leq_T Y$, where $Y$ is Martin-L\"of random with respect to some $Z\in\reals$, then $X$ is also Martin-L\"of random with respect to $Z$.  
\end{thm}

The phenomenon described in the ${XY\!Z}$ theorem can be seen as a kind of preservation of randomness:  if we map a random sequence $Y$ to a random sequence $X$, this map preserves relative randomness, in the sense that if an oracle $Z$ does not detect $Y$ as non-random, then it does not detect $X$ as non-random.

It immediately follows from this result that every Martin-L\"of random sequence Turing reducible to an $n$-random sequence for $n\geq 2$ is itself $n$-random.

\subsection{Chaitin's $\Omega$}  In \cite{c75}, Chaitin introduced his celebrated number $\Omega$, the halting probability, which is defined to be
\[
\Omega=\sum_{U(\sigma){\downarrow}}2^{-|\sigma|},
\]
where $U$ is a fixed universal, prefix-free machine.   That is, $\Omega$ is the probability that $U$ halts on some initial segment of an infinite input sequence.

Given that this definition is dependent on the choice of universal machine $U$, it is more accurate to define the family of $\Omega$-numbers $\Omega_U$ for every universal prefix-free $U$.  However, we will still refer to a fixed $\Omega$-number as $\Omega$, with the understanding that any property that we attribute to this $\Omega$-number will hold of all $\Omega$-numbers.

  Clearly $\Omega$ is left-c.e., since it is the limit of a computable sequence of partial sums determined by running $U$ for only finitely many steps.  More significantly, Chaitin proved that $\Omega$ is both Turing complete and Martin-L\"of random.  That is, $\Omega$ encodes all of the information of the halting problem, and yet the bits of $\Omega$ are arranged in such a way as to pass all Martin-L\"of tests.

Solovay, in unpublished work in the 1970s \cite{solovay75}, introduced a new reducibility for left-c.e.\ reals that is now known as \emph{Solovay reducibility}.  Given left-c.e.\ reals $\alpha$ and $\beta$, $\alpha$ is \emph{Solovay reducible} to $\beta$, written $\alpha\leq_S\beta$, if there is a $c\in\omega$ and a partial computable function $f:\mathbb{Q}\rightarrow\mathbb{Q}$ such that for every $q\in\mathbb{Q}$ such that $q<\beta$, $f(q)$ is defined, $f(q)<\alpha$, and $\alpha -f(q)<c(\beta-q)$.  That is, from any rational less than $\beta$ that is sufficiently close to $\beta$, we can compute a rational less than $\alpha$ that is sufficiently close to $\alpha$.  One consequence of this definition, the formal details of which we do not consider here, is that any computable sequence of rationals converging to $\beta$ from below can be effectively transformed into a computable sequence of rationals converging to $\alpha$ from below at the same rate (up to a multiplicative constant).

Solovay proved that for every left-c.e.\ real $\alpha$, $\alpha\leq_S\Omega$.  Calude, Hertling, Khoussainov, and Wang \cite{chkw01} later studied $\Omega$-like left-c.e.\ reals, where a left-c.e.\ real $\beta$ is \emph{$\Omega$-like} if for every left-c.e.\ real $\alpha$, $\alpha\leq_S\beta$.  The main result in \cite{chkw01} is that if $\beta$ is $\Omega$-like, then there is some universal, prefix-free machine $U$ such that $\beta=\Omega_U$.  Thus, the collection of $\Omega$-numbers and the collection of $\Omega$-like reals coincide.

The final piece in understanding the relationship between left-c.e.\ reals and Solovay reducibility was provided by Ku\v cera and Slaman, who proved the following:

\begin{thm}[Ku\v cera and Slaman \cite{ks01}]
If $\alpha$ is left-c.e.\ and Martin-L\"of random, then $\alpha$ is $\Omega$-like.
\end{thm}

Summing up, we have:

\begin{cor}  For a left-c.e.\ real $\alpha$, the following are equivalent:
\begin{itemize}
\item[(i)] $\alpha=\Omega_U$ for some universal, prefix-free machine $U$.
\item[(ii)] $\beta\leq_S\alpha$ for every left-c.e.\ real $\beta$, that is, $\alpha$ is $\Omega$-like.
\item[(iii)] $\alpha$ is Martin-L\"of random.
\end{itemize}
\end{cor}

For a survey on more recent work on $\Omega$, see Barmpalias's article in this volume.

\subsection{Randomness-theoretic reducibilities}

In addition to Solovay reducibility, a number of other reducibilities defined in terms of algorithmic randomness were introduced in the period under consideration.  We briefly survey four of them.  

First, two reducibilities originally introduced by Solovay \cite{solovay75}, given in terms of the comparison of the initial segment complexities of the sequences, are defined as follows: for some constant $c$, for every $n$
\[
A\leq_K B \Leftrightarrow K(A\restr n)\leq K(B\restr n)+c
\]
and
\[
A\leq_C B \Leftrightarrow C(A\restr n)\leq C(B\restr n)+c.
\]
Thus, for each of these reducibilities, the more random a sequence is, the higher it will be with respect to the orderings $\leq_K$ and $\leq_C$.
The resulting degree structures defined in terms of these reducibilities, the $K$-degrees and the $C$-degrees, respectively, have been studied in, for instance, \cite{dhns03}, \cite{ms07}, \cite{milleryu08}, and \cite{miller09}.  

Some key results on the structures of both the $K$-degrees and the $C$-degrees follow:
\begin{itemize}
\item For $\alpha,\beta\in\reals$, $\alpha\leq_C\beta$ implies $\alpha\leq_T\beta$, but $\alpha\leq_K\beta$ does not imply $\alpha\leq_T\beta$ (the former is attributed to Stephan in \cite[Theorem 9.7.1]{dhbook}).
\item Both the $C$-degrees and the $K$-degrees of left-c.e.\ reals form upper semilattices with the join operation given by addition \cite{dhns03}. 
\item There is an uncountable $K$-degree (attributed to J.\ Miller in \cite{b13}).
\item There is a minimal pair in the $K$-degrees \cite{cm06}, \cite{ms07}, \cite{bv11}.
\item There is a minimal $C$-degree \cite{ms07}.
\item There is a pair of $K$-degrees with no upper bound \cite{milleryu08}.
\end{itemize}

Two additional reducibilities, which are comparisons of the sequences' strengths as oracles, are:
\[
A\leq_{LR}B \Leftrightarrow \MLR^B\!\subseteq\MLR^A
\]
where, for $X\in\reals$, $\MLR^X$ stands for the collection of Martin-L\"of random sequences relative to $X$, and
\[
A\leq_{LK}B \Leftrightarrow \exists c \forall \sigma ( K^B(\sigma)\leq K^A(\sigma)+c).
\]
Both $\leq_{LR}$ and $\leq_{LK}$ were introduced by Nies in \cite{nies05}. Informally, the idea behind these reducibilities is that if $A$ is sufficiently powerful to detect that some object (a sequence $A\in\reals$ in the case of $\leq_{LR}$ and a string $\sigma\in\strings$ in the case of $\leq_{LK}$) is not random, then $B$ also detects that this object is not random.

Although $\leq_{LR}$ and $\leq_{LK}$ concern the relative power of oracles for determining the randomness of different kinds of objects (sequences and strings), remarkably, these reducibilities coincide.

\begin{thm}[Kjos-Hanssen, J.\ Miller, and Solomon \cite{kms12}]
For $A,B\in\reals$, 
\[
A\leq_{LR}B \text{ if and only if } A\leq_{LK}B.
\]
\end{thm}

Other significant results on the $LR$-degrees are:

\begin{itemize}
\item The set $\{A\st A\leq_{LR}\emptyset'\}$ is uncountable \cite{blsoskova08}.
\item The set $\{A\st A\leq_{LR} B\}$ is countable if and only if $\Omega$ is Martin-L\"of random with respect to $B$ \cite{bl11}.
\item There is some $A<_T\emptyset'$ such that $\emptyset'\leq_{LR}A$ \cite{cgm06} (see also \cite[Theorem 6.7]{simp07}).
\end{itemize}

These reducibilities allow us, as we will see in the next section, to calibrate sequences in terms of how far they are from being random.  For a thorough discussion of the above-discussed reducibilities, as well as others, see \cite{b13}.

\subsection{Other randomness notions and lowness for randomness}\label{subsec-lowness}

As the Turing degrees of \ml random sequences became better understood, there were newfound emphases on sequences that have properties that are antithetical to properties that \ml random sequences have and on understanding other randomness notions more completely. The first property is based on the reducibility $\leq_K$ discussed above and states that a sequence is far from random if it has low prefix-free Kolmogorov complexity:

\begin{definition}
A sequence $A$ is \emph{$K$-trivial} if there is a constant $c$ such that for all $n\nat$,
$$K(A\restr n) \leq K(0^n) + c.$$
\end{definition}

\noindent In other words, the $K$-trivials are the elements of the least $K$-degree.  

The second property states that a sequence is far from random if it is not capable of signficantly reducing the Kolmogorov complexity of any string, making it a lowness notion:

\begin{definition}(Muchnik, unpublished)
A sequence $A$ is \emph{low for $K$} if there is a constant $c$ such that for all strings $\sigma$,
$$K(\sigma)\leq K^A(\sigma) + c.$$
\end{definition}

We observe that $A$ is low for $K$ exactly when $A\leq_{LK} \emptyset$. Lowness was also defined in two other ways that are more generalizable to other randomness notions.

\begin{definition}
Let $\mathcal{R}$ be a randomness notion; in the second part of this definition, we assume that $\mathcal{R}$ has a test definition.
\begin{enumerate}
\item A sequence $A$ is \emph{low for $\mathcal{R}$-randomness} if every $\mathcal{R}$-random set is still $\mathcal{R}$-random relative to $A$ (in other words, if $A\leq_{LR} \emptyset$). 
\item A sequence $A$ is \emph{low for $\mathcal{R}$-tests} if, for every $\mathcal{R}$-test relative to $A$ $\langle V_i^A\rangle_{i\in\omega}$, there is an unrelativized $\mathcal{R}$-test $\langle U_i\rangle$ that covers it; e.g., $\bigcap_{i\in\omega} V_i^A \subseteq \bigcap_{i\in\omega} U_i$.
\end{enumerate}
\end{definition}

We note that, since there is a universal prefix-free machine and a universal \ml test, we automatically have that lowness for \ml randomness is equivalent to lowness for \ml tests.

The last standard way in which a sequence is considered weak with respect to a randomness notion is that of being a base. This definition of Ku\v{c}era, which first appeared in \cite{k93}, is also easily stated in a general way. It is based on Sacks' Theorem \cite{sacksdou}, which, as mentioned in Section \ref{subsec:demkuc}, states that if $A$ is not computable, then $\mu(\{ X\st X\geq_T A\})=0$. Here, instead, Ku\v{c}era states that a set should be far from random if something that computes it can be random relative to it.

\begin{definition}
Let $\mathcal{R}$ be a randomness notion. A sequence $A$ is a \emph{base for $\mathcal{R}$-randomness} if there is some $X\geq_T A$ such that $X$ is $\mathcal{R}$-random relative to $A$.
\end{definition}

Most of these lowness notions were proven to coincide by Nies in the case of \ml randomness in \cite{nies05}, when he proved that a sequence is low for \ml randomness exactly when it is low for $K$. In the same paper, Nies proved that $K$-triviality and lowness for \ml are equivalent; the proof that $K$-triviality implies lowness for $K$ is joint with Hirschfeldt. Ku\v{c}era showed in \cite{k93} that every sequence that is low for \ml randomness is a base for \ml randomness; the characterization was completed by Hirschfeldt, Nies, and Stephan in \cite{hns07}. We summarize these results here:

\begin{thm}[\cite{hns07,k93,nies05}]
The following are equivalent for a sequence $A$.
\begin{enumerate}
\item $A$ is low for \ml randomness.
\item $A$ is low for \ml tests.
\item $A$ is low for $K$.
\item $A$ is $K$-trivial.
\item $A$ is a base for \ml randomness.
\end{enumerate}
\end{thm}

The Turing degrees of these sequences have also been studied in great depth, beginning with Solovay's construction of a noncomputable $K$-trivial set in the 1970s \cite{solovay75}. Chaitin \cite{c77} and then Zambella \cite{za90} had shown that all $K$-trivial sequences are $\Delta^0_2$; in fact, Zambella used techniques based on Solovay's to show that there is a noncomputable $K$-trivial \ce set. In 1999, Ku\v{c}era and Terwijn proved that there is a noncomputable \ce set that is low for \ml randomness \cite{kt99}, and Muchnik proved that there is a noncomputable \ce set that is low for $K$ (unpublished). Later, Nies showed that when we consider only the \ce $K$-trivial Turing degrees, we get a $\Sigma^0_3$ ideal in the Turing degrees and that the $K$-trivial degrees form a Turing ideal that is generated by the \ce $K$-trivial degrees, in other words, that the smallest ideal containing the \ce $K$-trivial degrees consists of precisely the $K$-trivial degrees \cite{nies05}. In fact, he also showed that the $K$-trivial sets form a proper subclass of the superlow sets \cite{nies05}.

Of course, these definitions can be discussed in the context of other randomness notions. Work on Schnorr lowness was contemporary with the earliest work on lowness for \ml randomness. In \cite{ak00}, Ambos-Spies and Ku\v{c}era asked whether lowness for Schnorr randomness and lowness for Schnorr tests coincided; since there is no universal Schnorr test, the answer is not obviously yes. The first step in solving this question was taken by Terwijn and Zambella in \cite{tz01}, who characterized the Turing degrees that are low for Schnorr tests in terms of a concept they defined, computable traceability. This concept requires us to define a standard way to list the finite sets computably. While there are several ways to fix a computable list of the finite sets, here we will say that the $n^{th}$ canonical finite set $D_n$ contains precisely the natural numbers $i$ such that there is a 1 in the $2^i$'s place in the binary representation of $n$.

\begin{definition}
A Turing degree $\degd$ is \emph{computably traceable} if there is some order function $p$ such that for each function $f\leq_T \degd$, there is a computable function $r$ such that for all $n$,
\begin{enumerate}
\item $f(n)\in D_{r(n)}$ and
\item $|D_{r(n)}|\leq p(n)$,
\end{enumerate}
where $D_n$ is the $n^{th}$ canonical finite set. 

If we replace $D_{r(n)}$ with $W_{r(n)}$, we have a \emph{\ce traceable} degree.
\end{definition}

Now we can state Terwijn and Zambella's result:

\begin{thm}[Terwijn and Zambella, \cite{tz01}]
A Turing degree is low for Schnorr tests if and only if it is computably traceable.
\end{thm}

The corresponding result on lowness for Schnorr randomness was proven via contributions from a large number of people. We note here that every sequence that is low for Schnorr tests must be low for Schnorr randomness as well. Bedregal and Nies showed that any degree that is low for Schnorr randomness or for computable randomness must be hyperimmune-free \cite{bn03}; Kjos-Hanssen, Nies, and Stephan showed that a sequence $A$ is  \ce traceable exactly when every Schnorr null set relative to $A$ is \ml null as well and then that any hyperimmune-free \ce traceable degree must be computably traceable \cite{kns05}. We can combine these results and see the following:

\begin{thm}[Kjos-Hanssen, Nies, and Stephan \cite{kns05}]
The following are equivalent for a sequence $A$.
\begin{enumerate}
\item $A$ is low for Schnorr randomness.
\item $A$ is low for Schnorr tests.
\item $A$ is computably traceable.
\end{enumerate}
\end{thm}

Before Schnorr triviality could be considered, it was necessary to characterize Schnorr randomness in terms of Kolmogorov complexity. This required a new type of machine that was first defined by Downey and Griffiths in \cite{dg04}.

\begin{definition}
A \emph{computable measure machine} is a prefix-free machine $M$ such that the measure of $[\dom(M)]$ is a computable real.
\end{definition}

Downey and Griffiths then used computable measure machines to characterize Schnorr randomness; the characterization is directly parallel to that for \ml randomness. However, we note that it is inherently more complicated since there is no universal computable measure machine.

\begin{thm}[Downey and Griffiths \cite{dg04}]
A sequence $A$ is Schnorr random if and only if for every computable measure machine $M$, there is a constant $c$ such that for all $n$,
$$K_M(A\restr n)\geq n-c.$$
\end{thm}

Now we can discuss Schnorr triviality. In \cite{dg04}, Downey and Griffiths also defined a reducibility parallel to $\leq_K$ that is based on computable measure machines:

\begin{definition}
Let $A$ and $B$ be sequences. $A$ is \emph{Schnorr reducible to} $B$ (written $A\leq_{Sch} B$) if for every computable measure machine $M$, there is a computable measure machine $N$ and a constant $c$ such that for all $n$,
$$K_N(A\restr n)\leq K_M(B\restr n) + c.$$
\end{definition}

In a parallel with $K$-triviality, Downey and Griffiths defined a sequence $A$ to be \emph{Schnorr trivial} if $A\leq_{Sch} 0^\omega$. Schnorr triviality and lowness for Schnorr were soon seen not to coincide: almost immediately after Schnorr triviality was defined, Downey, Griffiths, and LaForte proved that there is a Turing complete Schnorr trivial sequence \cite{dgl04}. Later, Franklin proved that, in fact, every high degree contains a Schnorr trivial sequence \cite{f08-1}. This lets us see immediately that the Schnorr trivial sequences and the Schnorr low sequences do not coincide; however, Franklin showed that the hyperimmune-free Schnorr trivial degrees are precisely the degrees that are low for Schnorr \cite{f08-2}.

We may also remark that the Schnorr trivial degrees do not form an ideal in the Turing degrees: Downey, Griffiths, and LaForte also showed that there is a \ce degree that contains no Schnorr trivial (or $K$-trivial) sets \cite{dgl04}, so the Schnorr trivial Turing degrees are not closed downward. However, they did show that their truth-table degrees are closed downward, and since Franklin and Stephan showed that Schnorr trivial sequences are closed under join \cite{fsStrivtt}, we can at least say that they form an ideal in the truth-table degrees. For a more detailed discussion of the necessity of choosing the proper reducibility in randomness, see Franklin's survey.

The bases for Schnorr randomness form a different class yet again: Franklin, Stephan, and Yu proved that a sequence is a base for Schnorr randomness exactly when it cannot compute $\emptyset'$ \cite{fsy-bases}.

Just as Schnorr randomness was newly characterized in terms of Kolmogorov complexity in this period, computable randomness was newly characterized in terms of tests in two different ways. Downey, Griffiths, and LaForte developed the concept of a computably graded test \cite{dgl04}:

\begin{definition}
A \emph{computably graded test} is a pair $(\langle V_i\rangle_{i\in\omega}, f)$ where $\langle V_i\rangle_{i\in\omega}$ is a \ml test and $f:\strings\times \omega \rightarrow \R$ is a computable function such that the following three conditions hold for all $n\nat$, all $\sigma\in\strings$, and any finite prefix-free set of strings $\{\sigma_i\}$ such that $\bigcup[\sigma_i]\subseteq[\sigma]$:
\begin{enumerate}
\item $\mu(V_n \cap [\sigma])\leq f(\sigma,n)$,
\item $\sum f(\sigma_i,n) \leq 2^{-n}$, and
\item $\sum f(\sigma_i,n) \leq f(\sigma, n)$.
\end{enumerate}
A sequence passes a computably graded test if it passes the \ml test component.
\end{definition}

Merkle, Mihailovi\'{c}, and Slaman developed another test notion they used to characterize computable randomness, a bounded \ml test \cite{mms06}. They began by defining a \emph{computable rational probability distribution}: a computable function $\nu:\strings\rightarrow \Q$ such that $\nu(\langle\rangle)=1$ and $\nu(\sigma) = \nu(\sigma0)+\nu(\sigma1)$ for every $\sigma\in\strings$.

\begin{definition}
A \emph{bounded \ml test} is a \ml test $\langle V_i\rangle_{i\nat}$ such that there is a computable rational probability distribution such that for all $n$ and $\sigma$,
$$\mu(V_n\cap [\sigma])\leq 2^{-n}\nu(\sigma).$$ 
\end{definition}

\begin{thm}[Downey, Griffiths, and LaForte \cite{dgl04}, Merkle, Mihailovi\'{c}, and Slaman \cite{mms06}]
The following are equivalent for a sequence $A$.
\begin{enumerate}
\item $A$ is computably random.
\item $A$ passes all computably graded tests.
\item $A$ passes all bounded \ml tests.
\end{enumerate}
\end{thm}

There is also a machine characterization of computable randomness due to Mihailovi\'{c}; while his work is unpublished, the following definition and theorem appear in Section 7.1.5 of \cite{dhbook}.

\begin{definition}
A \emph{bounded machine} is a prefix-free machine such that there is a computable rational probability distribution $\nu$ such that for all $n$ and all $\sigma$,
$$\mu([\{\sigma \st K_M(\sigma)\leq |\sigma|-n\}]) \leq 2^{-n}\nu(\sigma).$$
\end{definition}

\begin{thm}
A sequence $A$ is computably random if and only if for every bounded machine $M$, there is a constant $c$ such that for all $n$, 
$$K_M(A\restr n)\geq n-c.$$
\end{thm}

Now we can turn our attention to weakness with respect to computable randomness. The first concept to be addressed was that of lowness, and Nies proved that, once again, we get a different class of sequences.

\begin{thm}[Nies, \cite{nies05}]
The sets that are low for computable randomness are precisely the computable sets.
\end{thm}

Lowness for tests for computable randomness has not been studied, nor has any notion of ``computable triviality." However, there is a partial characterization of the bases for computable randomness. Hirschfeldt, Nies, and Stephan have shown that every $\Delta^0_2$ sequence of non-DNC degree is a base for computable randomness but that no sequence of PA degree is \cite{hns07}. 

Now we turn our attention to randomness notions that are either stronger than \ml randomness or weaker than Schnorr randomness. In Wang's 1996 dissertation, he developed a \ml test-like characterization and a martingale characterization for Kurtz randomness \cite{wang}; this was followed almost a decade later by a similar characterization in terms of Kolmogorov complexity by Downey, Griffiths and Reid \cite{dnwy06}. 

\begin{definition}\cite{wang}
A \emph{Kurtz null test} is a \ml test $\langle V_i\rangle_{i\in\omega}$ such that for some computable function $f:\omega\rightarrow (\strings)^{<\omega}$, $V_n=[f(n)]$ for all $n$.
\end{definition}

\begin{definition}\cite{dgr04}
A prefix-free machine $M$ is a \emph{computably layered machine} if there is a computable function $f:\omega\rightarrow (\strings)^{<\omega}$ such that the following three conditions hold:
\begin{enumerate}
\item $\bigcup_{i\in\omega} f(i) = \dom(M)$,
\item if $\sigma\in f(i+1)$, then there is some $\tau\in f(i)$ such that $M(\tau)\preceq M(\sigma)$, and
\item if $\sigma\in f(i)$, then $|M(\sigma)| = |\sigma| + i+ 1$.
\end{enumerate}
\end{definition}

The equivalence of (1) and (2) in the following theorem is inherent in Kautz's thesis \cite{kautz} but explicitly stated in \cite{wang}; the equivalence of (4) to the others appears in \cite{dgr04}.

\begin{thm}[Wang \cite{wang}, Downey, Griffiths, and Reid \cite{dgr04}]
Let $A$ be a sequence. The following are equivalent.
\begin{enumerate}
\item $A$ is Kurtz random.
\item $A$ passes all Kurtz null tests.
\item For every computable martingale $d$ and every order function $h$, there is some $n$ such that $d(A\restr n)\leq h(n)$.
\item For every computably layered machine $M$, there is a constant $c$ such that for all $n$,
$$K_M(A\restr n) \geq n-c.$$
\end{enumerate}
\end{thm}

A comment must be made about lowness for Kurtz tests. Downey, Griffiths, and Reid proved that the degrees that are low for Kurtz tests form a superset of the computably traceable degrees and a subset of the hyperimmune-free degrees \cite{dgr04}; the final characterization would come later \cite{sy06,gm09}.

Now we turn our attention to weak 2-randomness. While weak 2-randomness had been briefly studied in the early 1980s by Gaifman and Snir \cite{gs82} and even earlier by Solovay, who proved that no weakly 2-random sequence is $\Delta^0_2$ \cite{solovay75}, it was placed into the randomness hierarchies we know today by Kurtz \cite{kurtz}, who defined weak $n$-randomness for all $n\geq 2$.  Later, weak 2-randomness was considered by Kautz \cite{kautz} and Wang \cite{wang}, who independently developed another characterization of it using a test definition similar to that of \ml randomness. 

\begin{definition}
A \emph{generalized \ml test} is a sequence $\langle V_i\rangle_{i\in\omega}$ of uniformly \ce subsets of $\reals$ such that $\lim_{n\rightarrow\infty} \mu(V_n)=0$.
\end{definition}

\begin{thm}[Kautz \cite{kautz}, Wang \cite{wang}]
A sequence is weakly 2-random if it passes all generalized \ml tests.
\end{thm}

However, a further study of the properties of weakly random sequences was not conducted until 2006, when Downey, Nies, Weber, and Yu \cite{dnwy06} undertook a systematic study of them.  Their main result is the following:

\begin{thm}[Downey, Nies, Weber, and Yu \cite{dnwy06}]
The following are equivalent for a \ml random sequence $A$.
\begin{enumerate}
\item $A$ is weakly 2-random.
\item The degree of $A$ and $\emptyset'$ form a minimal pair.
\item $A$ does not compute any noncomputable \ce set.
\end{enumerate}
\end{thm}

This gives us a characterization of the subclass of \ml random sequences that do not compute a noncomputable \ce set in terms of a randomness notion just as difference randomness characterizes the subclass of \ml random sequences that do not compute a PA degree \cite{fng}.

Downey, Nies, Weber, and Yu also began the study of lowness for weak 2-randomness by proving that every sequence that is low for weak 2-randomness is $K$-trivial \cite{dnwy06}; Kjos-Hanssen, J.\ Miller, and Solomon proved the converse in \cite{kms12}.

Now that we have discussed the Turing degrees of all of the main randomness notions, we can make an observation about randomness in the hyperimmune-free Turing degrees: all the main notions of randomness that we have seen thus far coincide in this setting.

\begin{thm}[Nies, Stephan, and Terwijn \cite{nst05}, Yu (unpublished)]
If $A$ has hyperimmune-free degree, then $A$ is Kurtz random if and only if $A$ is weakly 2-random.
\end{thm}

One additional result, also due to Nies, Stephan, and Terwijn \cite{nst05}, and independently, Miller \cite{miller04}, is a characterization of 2-randomness in terms of plain Kolmogorov complexity.  As stated in Section \ref{subsec-ml}, one cannot provide a definition of randomness by requiring all initial segments to be incompressible with respect to plain Kolmogorov complexity.  Martin-L\"of observed that almost every sequence $A$ satisfies the condition that there is some $c$ such that $C(A\restr n)\geq n-c$ for infinitely many $n$; let us say that such a sequence $A$ is \emph{infinitely often $C$-incompressible}.  The following result identifies precisely the notion of randomness that corresponds to this condition:

\begin{thm}[Nies, Stephan, and Terwijn \cite{nst05}, Miller \cite{miller04}]
For $A\in\reals$, $A$ is 2-random if and only if there $A$ is infinitely often $C$-incompressible.
\end{thm}

Miller \cite{miller04} further characterizes 2-randomness in terms of the property of being infinitely often $K$-incompressible, which is satisfied by a sequence $A$ if there is some $c$ such that $K(A\restr n)\geq n+K(n)-c$ for infinitely many $n$.

In this period, a randomness notion that was introduced in the late 1990s in \cite{msu98}, and which is still not entirely understood, began to receive interest: Kolmogorov-Loveland randomness. This notion is defined using a martingale characterization, and the key feature of this characterization is that the bits of the sequence do not need to be bet on in order. This idea has its roots in the notion of Kolmogorov-Loveland stochasticity \cite{kol63,kol98,loveland66} and was studied by Muchnik, Semenov, and Uspensky \cite{msu98} and Merkle, J.\ Miller, Nies, Reimann, and Stephan \cite{mmnrs06}; as in Downey and Hirschfeldt \cite{dhbook}, we use the latter group's notation. 

Let us say that a \emph{finite assignment} is a sequence of elements $(r_i,a_i)$ from $\omega\times \{0,1\}$ such that the $r_i$'s are all distinct: the $r_i$'s are the places on which bets have been made, and the $a_i$'s are the bets made at those places, so this gives a ``history" of how much has been bet and at which places. If we then define a \emph{scan rule}, which is a function that determines the next place to bet on given this ``history," and a \emph{stake function} that tells us how much to bet there, we have a \emph{nonmonotonic betting strategy}. We can say that a nonmonotonic betting strategy succeeds on a sequence $A$ if the limsup of the amount one has after betting on $A$ following the strategy defined by the scan rule and stake function is infinite.

\begin{definition}\cite{msu98}
A sequence is \emph{Kolmogorov-Loveland random} if no partial computable nonmonotonic betting strategy succeeds on it.
\end{definition}

Merkle would later prove that one could get the same class of random sequences using total computable nonmonotonic betting strategies \cite{merkle03}, so it is clear that every Kolmogorov-Loveland random is computably random; Muchnik, Semenov, and Uspensky showed that every \ml random sequence is Kolmogorov-Loveland random \cite{msu98}. It has been conjectured that the Kolmogorov-Loveland random sequences are precisely the \ml random sequences. While quite some time and energy has been put into this question \cite{kl10}, it remains open. 

\subsection{Effective notions of dimension}

Another significant strand of research in algorithmic randomness that emerged in the early 2000s involved effective notions of dimension. Lutz first provided a definition of effective Hausdorff dimension in \cite{lutz03}, given in terms of certain betting strategies called $s$-\emph{gales}.  First, Lutz gave a betting characterization of classical Hausdorff dimension in terms of these strategies.  Derived from the definition of a martingale given above, an \emph{$s$-gale} is a function $d:\strings\rightarrow\mathbb{R}^{\geq 0}$ satisfying $d(\sigma)=2^{-s}(d(\sigma0)+d(\sigma1))$ for all $\sigma\in\strings$.  The set of all sequences on which an $s$-gale $d$ succeeds (that is, those sequences $A$ such that $d(A\restr n)$ is unbounded in $n$) is written $S[d]$ as it is for a standard martingale.  Then Lutz proved the following:
\begin{thm}[Lutz \cite{lutz03}]
For $A\subseteq\reals$, 
\[
\mathrm{dim}_\mathrm{H}(A)=\inf\{s\in\mathbb{Q}\st A\subseteq S[d] \text{ for some $s$-gale }d\}.
\]
\end{thm}
For additional characterizations of dimension in terms of $s$-gales, see \cite{lutz03}.  To obtain effective Hausdorff dimension, we simply restrict our attention to c.e.\ $s$-gales, that is, those $s$-gales that take on values that are uniformly computable approximable from below.  We thus define
\[
\mathrm{dim}(A)=\inf\{s\in\mathbb{Q}\st A\subseteq S[d] \text{ for some c.e.\ $s$-gale }d\}
\]
to be the effective Hausdorff dimension of $A\subseteq\reals$.  One can equivalently define the effective Hausdorff dimension in terms of covers; see, for instance, the treatment in \cite[Section 13.5]{dhbook} for details.

One key difference between classical Hausdorff dimension and effective Hausdorff dimension is that individual sequences can have positive effective dimension, while sequences always have classical dimension zero.  Perhaps the most useful characterization of effective Hausdorff dimension for individual sequences was provided by Mayordomo, who proved the following:

\begin{thm}[Mayordomo \cite{mayordomo2002}]
For $A\in\reals$,
\[
\mathrm{dim}(A)=\liminf_{n}\dfrac{K(A\restr n)}{n}=\liminf_n\dfrac{C(A\restr n)}{n}.
\]
\end{thm}
Another fundamental early result about effective Hausdorff dimension is that for every computable real number $r\in[0,1]$, there is some $A\in\reals$ such that $\mathrm{dim}(A)=r$, a fact established by Lutz in \cite{lutz03}.

A similar characterization of an effective version of the classical notion of packing dimension can be given in terms of c.e.\ $s$-gales and initial segment complexity \cite{ahlm07}.  For the sake of brevity, we only highlight the latter.  As shown by Athreya, Hitchcock, Lutz, and Mayordomo \cite{ahlm07}, the effective packing dimension of a sequence $A\in\reals$, written $\mathrm{Dim}(A)$, can be characterized as 
\[
\mathrm{Dim}(A)=\limsup_{n}\dfrac{K(A\restr n)}{n}=\limsup_n\dfrac{C(A\restr n)}{n}.
\]

A number of computability theoretic questions about the notions of effective dimension have been investigated since these definitions first appeared, particularly in the context of the Turing degrees.  The most significant work in this respect has been on the broken dimension problem:  Given $\alpha\in(0,1)$, is there some sequence $A$ such that $\mathrm{dim}(A)=\alpha$ and for every $B\leq_TA$, $\mathrm{dim}(B)\leq\alpha$?  This question was answered in the affirmative by J.\ Miller, who showed that such a sequence $A$ exists and that it can even be chosen to be $\Delta^0_2$ \cite{miller11}.  Another important result due to Staiger \cite{staig98} and, independently, Hitchcock \cite{h05}, the correspondence principle for effective dimension, establishes conditions under which effective Hausdorff dimension and classical Hausdorff dimension agree:  For any countable union $U$ of $\Pi^0_1$ classes, $\mathrm{dim}_\mathrm{H}(U)=\mathrm{dim}(U)$.  

Lastly, we mention several recent uses of effective dimension in establishing new results on classical Hausdorff dimension.  These applications make use of the so-called \emph{point-to-set principle} due to J.\ Lutz and N.\ Lutz \cite{lutz2018}, which allows one to calculate the classical Hausdorff dimension of a subset $E$ of Euclidean space in terms of the relativizations of the effective dimensions of members of $E$.  That is, for $E\subseteq\mathbb{R}^n$,
\[
\mathrm{dim}_H(E)=\min_{A\in\reals}\sup_{x\in E}\mathrm{dim}^A(x),
\]
where, for $A\in\reals$ and $x\in\mathbb{R}^n$, $\mathrm{dim}^A(x)$ is the effective dimension of $x$ relative to the oracle $A$ (which can be defined in terms of $K^A$, prefix-free Kolmogorov complexity relative to $A$).  Using the point-to-set principle, N.\ Lutz and Stull \cite{stlutz2017} provided a new lower bound for the classical dimension of generalized Furstenberg sets, and N.\ Lutz \cite{lutz2017} used the principle to extend an inequality bounding the dimension of certain intersections of sets, known to hold for Borel subsets of Euclidean space, to all subsets of Euclidean space.

For more results on effective notions of dimension, see \cite[Chapter 13]{dhbook}.

\section{Recent developments}

Since the publication of the two monographs due to Nies \cite{niesbook} and Downey and Hirschfeldt \cite{dhbook}, a majority of research in algorithmic randomness (at least outside of the setting of resource-bounded randomness) has involved notions of randomness that fall between Kurtz randomness and 2-randomness when measured in terms of strength.  One exception is UD-randomness, introduced by Avigad \cite{av13} and further developed by Calvert and Franklin \cite{cf15}.  Based on the concept of uniform distribution studied by Weyl \cite{weyl16}, UD-randomness is implied by Schnorr randomness but incomparable with Kurtz randomness, as shown by Avigad.  

In particular, there has been a great deal of interest in notions even stronger than Martin-L\"of randomness.
Demuth randomness and weak Demuth randomness have already been discussed.  As mentioned in Section \ref{subsec:randseqTdegs}, in 2011, Franklin and Ng introduced the notion of difference randomness and showed that the sequences that are difference random are precisely the Turing incomplete \ml random sequences \cite{fng}.  Combined with Theorem \ref{stephan-thm} due to Stephan, it follows that a Martin-L\"of random sequence is difference random if and only if it does not compute a completion of Peano arithmetic.  Since then, a multitude of notions, including balanced randomness \cite{fhmnn10}, Oberwolfach randomness \cite{bgknt16}, and density randomness \cite{mnz16}, has been studied. Some of them have been used to characterize subclasses of the \ml random sequences in terms of computational strength, while others have proven useful for researchers who wish to relate algorithmic randomness to computable analysis \cite{mnz16}.

One of the primary achievements of recent work on the interaction between algorithmic randomness and computable analysis has been to provide characterizations of various notions of randomness in terms of almost everywhere behavior from classical mathematics, a research project that began in the work of Demuth but was only recently independently rediscovered.  Given a theorem of the form ``For almost every $x$, $x$ has property $P$," we can often replace the property $P$ with a computably defined analogue $P^*$ and prove a theorem of the form ``$x$ is $\mathcal{R}$-random if and only if $x$ has the property $P^*$", where $\mathcal{R}$ is some notion of effective randomness. This has been done for theorems such as the Lebesgue differentiation theorem \cite{bmn16,mnz16} and Birkhoff's ergodic theorem and the Poincar\'{e} recurrence theorem \cite{v97,ghr11,fgmn12,bdhms12,ft-mp}.   Many of the most significant recent results involving computable randomness and Schnorr randomness have been achieved in this area.
 The survey by Hoyrup on layerwise computability, the survey by  Rute on the relationship between algorithmic randomness and analysis, and the survey by Towsner on ergodic theory and randomness in this volume address the ongoing research in these areas. 
 
 \section{Acknowledgments}
The authors would like to thank the referees who provided extraordinarily perceptive and useful comments on this survey.

\def\cprime{$'$}

\end{document}